\documentclass[11pt,reqno]{amsart}
\usepackage{amsfonts}
\usepackage{mathrsfs}
\usepackage{amssymb}
\input{epsf.sty}
\catcode `\@=11
\def\numberbysection{\@addtoreset{equation}{section}
         \renewcommand{\theequation}{\thesection.\arabic{equation}}}
\numberbysection
\def\subsubsection{\@startsection{subsubsection}{3}%
  \normalparindent{.5\linespacing\@plus.7\linespacing}{-.5em}%
  {\normalfont\bfseries}}
\newtheorem{thm}{Theorem}[section]
\newtheorem{lem}[thm]{Lemma}
\newtheorem{prop}[thm]{Proposition}
\newtheorem{cor}[thm]{Corollary}

\theoremstyle{definition}

\newtheorem{df}[thm]{Definition}
\newtheorem{qu}[thm]{Question}
\newtheorem{rmk}[thm]{Remark}
\newtheorem{nota}[thm]{Notation}

\def\nn{\nonumber}

\def\eps{\epsilon}

\def\e{{\rm Eu}}
\def\et{{\rm Eu}_t}

\def\Xm{X^m}

\def\bm{{\bf m}}
\def\F{{\mathcal F}}
\def\R{{\mathcal R}}
\def\T{{\mathcal T}}
\def\rk{{\mathrm rk}}
\def\Sn{{\mathbb S}_n}

\newcommand{\eu}[1]{\e_{#1}}
\newcommand{\eut}[1]{\e_{#1,t}}
\newcommand{\fe}{\mathfrak{Eu}}
\newcommand{\coeff}[1]{\text{Coeff of $#1$ in }}

\newcommand{\evalf}[1]{\text{eval}_{\F|#1}}
\def\g{{\bf g}}
\def\Form{\Upsilon}

\def\vr{{\mathrm vr}}
\newcommand{\FF}[2]{{\F}_{#2}(#1)}
\def\fS{{\mathfrak S}}
\def\L{\mathscr L}

\begin{document}

\title{A Note on the Two Approaches to Stringy Functors for Orbifolds}

\author
[Ralph M.\ Kaufmann]{Ralph M.\ Kaufmann}
\email{kaufmann@math.uconn.edu}

\address{University of Connecticut, Department of Mathematics,
Storrs, CT 06269}

\begin{abstract}
In this note, we reconcile two approaches that have been used to
construct stringy multiplications. The pushing forward after pulling
back that has been used to give a global stringy extension of the
functors $K_0,K^{top}, A^*,H^*$\cite{CR,FG,AGV,JKK2}, and the
pulling back after having pushed forward, which we have previously
used in our (re)-construction program for $G$--Frobenius algebras,
notably in considerations of singularities with symmetries and for
symmetric products. A similar approach was also used by
\cite{ChenHu} in their considerations of the Chen--Ruan product in a
deRham setting for Abelian orbifolds.

We show that the pull--push formalism has a solution by the
push--pull equations in two situations. The first is a deRham
formalism with Thom push--forward maps and the second is the setting
of cyclic twisted sectors, which was at the heart of the
(re)-construction program.

We go on to do formal calculations using fractional Euler classes
which allows us to formally treat all the stringy multiplications
mentioned above in the general setting. The upshot is the formal
trivialization of the co--cycles of the reconstruction program using
the presentation of the obstruction bundle of \cite{JKK2}.
\end{abstract}

\maketitle

\section*{Introduction}
For global quotients there is a by now standard approach to
constructing  stringy products via first pulling back and then
pushing forward \cite{CR,FG,AGV,JKK2}. We will call this
construction the pull--push, which stands for pull after pushing.
However, going back to \cite{wisc,orb} we have used a mechanism that
first pushes forward and then pulls back to construct and classify
$G$--Frobenius algebra structures. In the same spirit we call this
the push--pull
--- read push after pulling.
This approach has been very successful for singularities
\cite{orb,orbsing} and for special cases of the group, for instance
$G={\mathbb S}_n$, see \cite{sq}. The advantage of this approach is
that one is left with solving an algebraic co-cycle equation. In
many cases this cocycle is unique up to normalized discrete torsion
\cite{orb,hilb,sq,disc,orbsing}. In fact as we argued in
\cite{orb,sq} this mechanism must work if the twisted sectors are
cyclic modules over the untwisted sector. Surprisingly a similar
technique to ours was used in \cite{ChenHu} where the authors passed
to the deRham chains and used formal fractional Thom forms to study
the product.

The goal of this paper is to consolidate these various results:

First, in \S\ref{cyclicsection}, we show that in the case of cyclic
twisted sectors both approaches exist for all the geometric functors
considered in \cite{JKK2}. This means that our reconstruction
program of \cite{wisc,orb} (see also \cite{hilb} for a short
detailed version) has a solution. It actually then has at least
discrete torsion many \cite{disc}. The key in this situation is the
existence of sections of the pull--back maps which allow us to prove
the relevant theorems using only the projection formula.

Next, in \S\ref{deRhamsection}, we show that in the deRham setting,
all the elements of the previous study hold up to homotopy, that is
up to exact forms. Hence we can provide a rigorous setting using
Thom push--forwards and pull--backs for general global quotients.

In the general setting, see \S\ref{formalsection}, there are some
obstacles towards a rigorous calculus of pushing forward and then
pulling back. We can make a lot of headway using the excess
intersection formula. But then we have to deal with formal
fractional Euler classes. There are two types of these classes. The
first are formal Euler classes of negative bundles, to be precise
the $-1$ multiples of the normal bundles of the fixed point sets
considered to live in $K$--theory. We can make sense of these as
formally defining sections and in the situations above these
sections coincide with the ones we constructed. The second type of
formal class is that of the fractional Euler class of positive, but
fractional classes in rational $K$--theory. This type of class poses
less of a problem and can be treated by adjoining roots to the
various rings.

One main result of the formal and rigorous calculations is that in
the different situations the classes $S_{m}$ appearing in the
definition of the obstruction bundle produce a co--cycle in the
sense of \cite{wisc,orb} that is trivialized by them.

\section*{Acknowledgements}
This paper  came into existence due to the continued interest and
questions of Yongbin Ruan,  whom we thank very much. It also owes a
lot to the ``Workshop on Quantum Cohomology of Stacks'' at the IHP
in Feb.~2007, where these results were first formulated explicitly
and presented. We wish to thank the organizers for the wonderful
conference and stimulating atmosphere. Last but not least we also
wish to thank Takashi Kimura for valuable conversations.

\section*{Conventions}
Will use at least coefficients in $\mathbb Q$ if nothing else is stated. For some applications such as deRham forms we will use $\mathbb R$ coefficients. All statements remain valid when passing to $\mathbb C$.

\section{General setup}
We will work in the same setup as in the global part of \cite{JKK2}.
That is we simultaneously treat two flavors of geometry, algebraic
and differential. For the latter, we consider a stably almost
complex manifold $X$ with the action of a finite group $G$ such that
the stably almost complex bundle is $G$ equivariant. While for the
former $X$ is taken to be a smooth projective variety.

In both situations for $m\in G$ we denote the fixed point set of $m$
by $X^m$ and let
\begin{equation}
I(X)=\amalg_{m\in G} X^m
\end{equation}
be the inertia variety.

We let ${\mathcal F}$ be any of the functors $H^*,K_0,A^*,K^{\rm
top}$, that is cohomology, Grothendieck $K_0$, Chow ring or
topological $K$--theory with ${\mathbb Q}$ coefficients,  and define
\begin{equation}
\F_{stringy}(X,G):= \F(I(X))= \bigoplus_{m\in G} \F(X^m)
\end{equation}
{\em additively}.

We furthermore set
\begin{equation}
\eu{\F}(E)=\begin{cases} c_{top}(E)&\text{ if } \F=H^* \text { or }
A^*
\text{ and $E$ is a bundle}\\
 \lambda_{-1}(E^*)&\text{ if } \F=K \text { or } K^{top}\\
 \end{cases}
\end{equation}

Notice that on bundles $\e$ is multiplicative. For general $K$--theory elements we set
\begin{equation}
\eut{\F}(E)=\begin{cases} c_{t}(E)&\text{ if } \F=H^* \text { or } A^*\\
 \lambda_{t}(E^*)&\text{ if } \F=K \text { or } K^{top}
\end{cases}
\end{equation}

\begin{rmk}
Notice $\eut{\F}$ is {\em always} multiplicative and it is a power
series that starts with $1$ and hence is invertible in $\F(X)[[t]]$.
\end{rmk}

\begin{df}
\label{evaldef}
 For a positive element $E$, i.e.\ $E$ can be
represented by a bundle with rank $r=\rk(E)$, we have that
$\eu{\F}(E)=\eut{\F}(E)|_{t=-1}$ for $\F$ either $K_0$ or $K^{top}$
and $\eu{\F}(e)=\coeff {t^r}[\eut{\F}(E)]$ if $\F$ is $A^*$ or
$H^*$. To be able to deal with both situations, for $E,r$ as above,
we define
\begin{equation}
\label{evaldefeq}
\evalf{r}(\eut{\F}(E))=\begin{cases}\eut{\F}(E)|_{t=-1} &
\text{ if $\F$ is $K_0$ or $K^{top}$}\\
\coeff {t^r}[\eut{\F}(E)]& \text{ if $\F$ is $A^*$ or $H^*$}
\end{cases}
\end{equation}
we then have $\evalf{r}(\eut{\F}(E))=\eu{\F}(E)$
\end{df}

\begin{rmk}
Notice that for $\F$ as above and each subgroup $H\subset G$,
$\F(X^H)$ is an algebra. We will call the internal product
$\F(X^m)\otimes \F(X^H)\to \F(X^H)$ the na\"ive product. There is
however a ``stringy--product'' which preserves the $G$--grading. To
define it, we recall some definitions from \cite{JKK2}.
\end{rmk}

\subsection{The stringy product via pull--push}
For $m \in G$ we let $\Xm$ be the fixed point set of $m$ and for a
triple $\bm=(m_1,m_2,m_3)$ such that  $\prod m_i={\bf 1}$ (where
${\bf 1}$ is the identity of $G$) we let $X^{\bm}$ be the common
fixed point set, that is the set fixed under the subgroup generated
by them.

In this situation, recall the following definitions. Fix $m \in G$ let $r=ord(m)$ be its order.
Furthermore let $W_{m,k}$ be the sub--bundle of $TX|_{X^{\bm}}$ on which $m$ acts with character $\exp(2\pi i \frac{k}{r})$, then
\begin{equation}
S_m=\bigoplus_k \frac{k}{r} W_{m,k}
\end{equation}
Notice this formula is invariant under stabilization.

We also wish to point out that using the identification $X^m=X^{m^{-1}}$
\begin{equation}
\label{normaleq}
S_m\oplus(S_{m^{-1}})=N_{X^m/X}
\end{equation}
where for an embedding $X\to Y$ we will use the notation $N_{X/Y}$
for the normal bundle.

Recall from \cite{JKK2} that in such a situation there is a product
on $\F(X,G)$ which is given by

\begin{equation}
\label{proddefeq}
 v_{m_1}*v_{m_2}:=\check e_{m_3*}(e^*_1(v_{m_1})e_2^*(v_{m_2})\e(\R(\bm)))
\end{equation}

where the obstruction bundle ${\mathcal R}(\bm)$ can be defined by
\begin{equation}
{\mathcal R}(\bm)=  S_{m_1}\oplus S_{m_2}\oplus S_{m_3} \ominus
N_{X^{\bm}/X}
\end{equation}
and the $e_i:X^{m_i}\to X$ and $\check e_3:X^{m_3^{-1}}\to X$ are
the inclusions. Notice, that as it is written $\R(\bm)$ only has to
be an element of K-theory with rational coefficients, but is
actually indeed represented by a bundle \cite{JKK2}.

\begin{rmk}
The first appearance of a push--pull formula was given in \cite{CR}
in terms of a moduli space of maps. The product was for the $G$
invariants, that is for the $H^*$ of the inertia orbifold and is
known as Chen--Ruan cohomology. In \cite{FG} the obstruction bundle
was given using Galois covers establishing a product for $H^*$ on
the inertia variety level, i.e.\ a $G$--Frobenius algebra as defined
in \cite{wisc,orb}, which is commonly referred to as the
Fantechi--G\"ottsche ring. In \cite{JKK1}, we put this global
structure back into a moduli space setting and proved the trace
axiom. The multiplication on the Chow ring $A^*$ for the inertia
stack was defined in \cite{AGV}. The representation of the
obstruction bundle in terms of the $S_m$ and hence the passing to
the differentiable setting as well as the two flavors of $K$--theory
stem from \cite{JKK2}.
\end{rmk}

%put back \check.

The following is the key diagram:

\begin{equation}
\label{maindiagrameq}
\begin{matrix}
&X&\\
i_1\nearrow &\uparrow i_2 &\nwarrow \check \imath_3\\
X^{m_1}&X^{m_2}&X^{m_3^{-1}}\\
e_1\nwarrow &\uparrow e_2 &\nearrow \check e_3\\
&X^{\bm}&
\end{matrix}
\end{equation}

Here we used the notation of \cite{JKK2}, where $e_3:X^{\bm}\to
X^{m_3}$ and $i_3: X^{m_3}\to X$ are the inclusion, $\vee: I(X)\to
I(X)$ is the involution which sends the component $X^m$ to
$X^{m^{-1}}$ using the identity map and $\check \imath_3= i_3\circ
\vee$, $\check e_3=\vee\circ e_3$. This is short hand notation for
the general notation of the inclusion maps $i_m:X^m\to X$, $\check
\imath_m:=i_m\circ \vee=i_{m^{-1}}$.

\subsection{The $\F(X)$ module structure}
Notice that each $\F(X^m)$ is an $\F(X)$ module in two ways which
coincide. First via the na\"ive product and pull back, i.e.\ $a\cdot
v_m:=i_m^*(a)v_m$ and secondly via the stringy multiplication
$(a,v_m) \mapsto a*v_m$. Now using (\ref{normaleq}) it is
straightforward to check that
\begin{equation}
a\cdot v_m = i_m^*(a)v_m = a*v_m
\end{equation}

\section{Pull--push: the cyclic case}
\label{cyclicsection}
 The way the product is defined in
(\ref{proddefeq}) is via first pulling back and then pushing forward
using the maps $e_k$. The aim of this section is to establish
rigorous arguments, that one can also first push--forward and then
pull back while using the maps $i_k$. This can be done rigorously
using sections and the projection formula. We apply this technique
in the current paragraph which treats the cyclic case and in
\S\ref{deRhamsection} which is devoted to the deRham setting.

\subsection{Sections}
We can realize the (re)--construction program of
\cite{wisc,orb,hilb} in two different situations. First, for any
functor $\F$ as above provided there are sections to the pull--back
maps $i_k^*$ and secondly in a deRham setting, where these sections
exist on the level of forms.

\begin{df}
We say that $\F$ admits sections for $(X,G)$
if for every map $i_m:X^m\to X$ there are sections $i_{ms}:\F(X^m)\to\F(X)$ of the pull--back maps
$i_{m}^*:\F(X)\to\F(X)$, that is $i^*_m \circ i_{ms}=id:\F(X^m)\to \F^(X^m)$
\end{df}

Examples are for instance given by symmetric products $(X^{\times n},\Sn)$, see \cite{orb,hilb} or
manifolds whose fixed loci are empty or points.

\begin{lem}
 If $\F$ admits sections for $(X,G)$, then $\F(X^m)$ is a cyclic $\F(X)$ module,
where the module structure is given by $a\cdot v_m:=i^*(a)v_m$. A
cyclic generator is $1_m$ which is the identity element of the
algebra $\F(X^m)$ endowed with the na\"ive product.

\end{lem}
\begin{proof}
$v_m=i^*_m(i_{ms}(v_m))= i_{ms}(v_m)\cdot 1_m$
\end{proof}

\begin{rmk}
We have

\begin{equation}
\label{ismulteq}
i^*_m(i_{ms}(a)i_{ms}(b))=i^*_m(i_{ms}(a))i^*_m(i_{ms}(b))=ab=
i^*_m(i_{ms}(ab))
\end{equation}
\end{rmk}

\subsection{A rigorous calculation using sections}

%In the case we have sections, the following calculation holds true.

\begin{prop}
\label{mainprop}
If there are sections $i_{js}$ of $i^*_j$ then the following
equation holds
\begin{eqnarray}
\label{sectionequation} && \coeff{t^r}\left\{\check
\imath^*_3[i_{1s}( v_{m_1}) i_{2s}
(v_{m_2})\gamma_{m_1,m_2}(t)]\right\}\nn\\
&=&\coeff{t^r}\left\{ \check \imath^*_3[i_{1s}(
v_{m_1}\et(S_{m_1}))i_{2s} (v_{m_2}\et (S_{m_2}))i_{3s}(\et (
S_{m_3})e_{3*}( \et(\ominus
N_{X^{\bm}/X})))]\right\}\nn\\
&=&v_{m_1}*v_{m_2}
\end{eqnarray}
where the product $*$ is the product defined in (\ref{proddefeq})
and $r=\rk(\R(\bm))$ and
\begin{eqnarray}
\label{cocycleeq} \gamma_{m_1,m_2}(t)&=&i_{1s}(\et(S_{m_1}))
i_{2s}(\et (S_{m_2}))i_{3s}(\et ( S_{m_3}) e_{3*}( \et(\ominus N_{X^{\bm}/X})))\nn\\
&=&i_{1s}(\et(S_{m_1}))i_{2s}(\et (S_{m_2}))\check \imath_{3s}(\et
(\ominus S_{m_3^{-1}})
\check e_{3*}( \et(\ominus N_{X^{\bm}/X^{m_3}})))\nn\\
\end{eqnarray}
\end{prop}

\begin{proof}
Using the projection formula, the defining equation for the sections
$i_j^*\circ i_{js}=id$, and the fact that $e_k\circ i_k=j=\check
e_k\circ \check \imath_k$
\begin{eqnarray}
&&\check e_{3*}[e_1^*(v_{m_1})e_2^*(v_{m_2}) \et(S_{m_1}|_{X^{\bm}}
\oplus S_{m_2}|_{X^{\bm}} \oplus S_{m_3}|_{X^{\bm}}\ominus
N_{X^{\bm}/X})]\nn\\
&=& \check e_{3*}[e_1^*(i^*_1(i_{1s}(
v_{m_1}\et(S_{m_1}))))e_2^*(i_2^*(i_{2s} (v_{m_2}\et
(S_{m_2}))))e_3^*(i^*_3(i_{3s}(\et ( S_{m_3})))) \et(\ominus
N_{X^{\bm}/X})]\nn\\
&=& \check e_{3*}[\check e_3^*(\check \imath^*_3(i_{1s}(
v_{m_1}\et(S_{m_1}))))\check e_3^*(\check \imath_3^*(i_{2s}
(v_{m_2}\et (S_{m_2}))))\check e_3^*(\check \imath^*_3(i_{3s}(\et (
S_{m_3})))) \et(\ominus
N_{X^{\bm}/X})]\nn\\
&=& \check \imath^*_3[i_{1s}( v_{m_1}\et(S_{m_1}))i_{2s} (v_{m_2}\et
(S_{m_2}))i_{3s}(\et ( S_{m_3})) \check \imath_{3s}(\check e_{3*}(
\et(\ominus
N_{X^{\bm}/X})))]\nn\\
&=& \check \imath^*_3[i_{1s}( v_{m_1}\et(S_{m_1}))i_{2s} (v_{m_2}\et
(S_{m_2}))i_{3s}((\et ( S_{m_3})) e_{3*}( \et(\ominus
N_{X^{\bm}}/X)))]\nn\\
\end{eqnarray}
So that taking the coefficient of $t^r$ with $r=\rk(\R(\bm))$ we
obtain the second claimed equality.
For the first equality we can use the fact (\ref {ismulteq})
\begin{eqnarray}
&&\check e_{3*}[e_1^*(v_{m_1})e_2^*(v_{m_2}) \et(S_{m_1}|_{X^{\bm}}
\oplus S_{m_2}|_{X^{\bm}} \oplus S_{m_3}|_{X^{\bm}}\ominus
N_{X^{\bm}/X})]\nn\\
&=& \check e_{3*}[e_1^*(i^*_1(i_{1s}(
v_{m_1}\et(S_{m_1}))))e_2^*(i_2^*(i_{2s} (v_{m_2}\et
(S_{m_2}))))e_3^*(i^*_3(i_{3s}(\et ( S_{m_3})))) \et(\ominus
N_{X^{\bm}/X})]\nn\\
&=& \check e_{3*}[e_1^*(i^*_1(i_{1s}(
v_{m_1})i_{1s}(\et(S_{m_1}))))e_2^*(i_2^*(i_{2s} (v_{m_2})i_{2s}
(\et (S_{m_2}))))\nn\\
&&\quad e_3^*(i^*_3(i_{3s}(\et ( S_{m_3})))) \et(\ominus
N_{X^{\bm}/X})]\nn\\
\end{eqnarray}
and proceed as above. Finally, for (\ref{cocycleeq}), we notice that
$N_{X^{\bm}/X}=N_{X^{\bm}/X^{m_3}}\oplus N_{X^{m_3}/X}|_{X^{\bm}}$
and use (\ref{normaleq}).
\end{proof}

%not quite true, need to set $t=1$

\begin{thm}
 Let $\F\in \{A^*,H^*,K_0,K^*_{\rm top}\}$ and
 $(X,G)$  in the appropriate category which admits sections for $\F$ then
 the equation (\ref{sectionequation}) solves the re--construction program of \cite{orb}
 with the co--cycles $\gamma_{m_1,m_2} :=\coeff{t^r}(\gamma_{m_1,m_2}(t))$.
\end{thm}

\begin{proof}
In this setting the calculation of the Proposition \ref{mainprop}
applies, which also shows, {\it a forteriori} that the formulas are
independent of the choice of lift and that the $\gamma_{m_1,m_2}
:=\coeff{t^r}(\gamma_{m_1,m_2}(t))$ are indeed co--cycles and
section independent co--cycles in the sense of \cite{orb}.
\end{proof}

\begin{rmk}
As we show in \S\ref{formalsection} below, these co--cycles are
formally trivial.
\end{rmk}

\subsection{Symmetric Product}
In particular the theorem above applies to  symmetric products and
gives a new way to show the existence of the unique co--cycles in
this situation constructed in \cite{sq}.

\section{The Chain level: a rigorous calculation using deRham Chains}
\label{deRhamsection} Although it is not true in general that the
pull back $e^*_i$ is surjective on cohomology or by the usual Chern
isomorphism on K--theory, on the level of deRham chains this is
true. Notice that in the proof of Proposition \ref{mainprop}, we
only used the projection formula, the defining equation for the
sections and the fact that the pull--back is an algebra
homomorphism.

\begin{nota}
In this section, we fix coefficients to be ${\mathbb R}$ and we
denote by $\Omega^n(X)$ the $n$--forms on $X$. Likewise for a bundle
$E\to B$ with compact base we denote $\Omega^n_{cv}(E)$ the $n$
forms on $E$ with compact vertical support and let $H^*_{cv}(E)$ the
corresponding cohomology with compact vertical support.
\end{nota}

\subsection{DeRham chains and Thom push--forwards} In this
section, we will use deRham chains and the Thom construction
\cite{BottTu}. The advantage is that every form on every $X^m$ is a
``pull--back'' from a tubular neighborhood.

We recall the salient features adapted to our situation from
\cite{BottTu}. Let $i:X\to Y$ be an embedding, then  there is a
tubular neighborhood $Tub(N_{X/Y})$ of the zero section of the
normal bundle $N_{X/Y}$ which is contained in $Y$. We let
$j:Tub(N_{X/Y})\to X$ be the inclusion.

Now the Thom isomorphism $\T: H^*(X)\to H^{*+codim(
X/Y)}_{cv}(N_{X/Y})$ can be realized on the level of forms via
capping with a  Thom form $\Theta$: $\T(\omega)=\pi^*(\omega)\wedge
\Theta$. The Thom map is inverse to the integration along the fibre
$\pi_*$ and hence $\pi_*(\Theta)=1$. In fact, the class of this form
is the unique class whose vertical restriction is a generator and
whose integral along the fiber is $1$. For any given tubular
neighborhood $Tub(N_{X/Y})$ of the zero section of the normal bundle
one can find a form representative $\Theta$ such that the
$supp(\Theta) \subset Tub(N_{X/Y})$.

\subsection{Push--forward}
In this situation the Thom push-forward $i_*: H^*(X)\to H^*(Y)$ is
given by $\T$ followed by the extension by zero $j_*$. These maps
are actually defined on the form level.  That is we choose  $\Theta$
to have support strictly inside the tube, and hence the extension by
zero outside the tube is well defined for the forms in the image of
the Thom map.
\begin{equation}
i_{*}(\omega):=j_*(T(\omega))=j_*(\pi^*(\omega)\wedge \Theta)
\end{equation}

Notice that for two consecutive embeddings
$X\stackrel{e}{\rightarrow}Y\stackrel{i}{\rightarrow}Z$, on
cohomology we have $e_*\circ i_*= (e\circ i)_*:H^*(X)\to H^*(Z)$. On
the level of forms depending on the choice of representatives of the
Thom form either the identity holds on the nose, since the Thom
classes are multiplicative \cite{BottTu} or they differ by an exact
form $e_*\circ i_*(\omega)= (e\circ i)_*+d\tau$.

\subsection{The projection formula on the level of forms}
The following proposition follows from standard facts \cite{BottTu}
\begin{prop}[Projection formula for forms]

With $i:X\to Y$ and embedding and $i_*$ defined as above, for any
form $\omega\in \Omega^*(X)$ and any closed form $\phi\in
\Omega^*(Y)$ there is an exact form $d\tau\in \Omega^*(Y)$ such that
\begin{equation}
i_*(i^*(\omega)\wedge \phi)= \omega \wedge i_*(\phi)+d\tau
\end{equation}
\end{prop}

\begin{proof}

 Denote the zero
section by $z:X\to N_{X/Y}$ and projection map of the normal
bundle by $\pi:N_{X/Y}\to X$, then $i=j\circ z$.
\begin{equation}
\begin{matrix}
X &{\stackrel{\pi|Tub}{\leftarrow}}\atop
{\stackrel{z}{\rightarrow}}&Tub(N_{X/Y})&
\stackrel{j}{\rightarrow}&Y
\end{matrix}
\end{equation}
Since $\pi$ is a deformation retraction,  $\pi^*$ and $z^*$ are
chain homotopic \cite{BottTu}, hence $\pi^*\circ z^* (\omega)=
\omega + d\tau$. We can now calculate
\begin{eqnarray}
i_*(i^*(\omega)\wedge \phi)&=&j_* (\pi^*(i^*(\omega)\wedge\phi)\wedge\Theta)\nn\\
&=&j_* (\pi^*(z^*(j^*(\omega))\wedge\pi^*(\phi)\wedge\Theta))\nn\\
&=&j_* ((j^*(\omega)+d\tau)\wedge\pi^*(\phi)\wedge\Theta)\nn\\
&=&\omega\wedge j_*(\pi^*(\phi)\wedge\Theta)+j_*(d\tau \wedge\pi^*(\phi)\wedge\Theta)\nn\\
&=&\omega\wedge i_*(\phi) +dj_*(\tau\wedge \pi^*(\phi)\wedge\Theta)
\end{eqnarray}
where the penultimate question holds true, since $\Theta$ has
support inside $Tub(N_{X/Y})$ and the last equation holds true since
$d$ commutes with the extension by zero and pull--back.
\end{proof}

\subsection{Sections}

\begin{figure}
\epsfxsize = \textwidth
\epsfbox{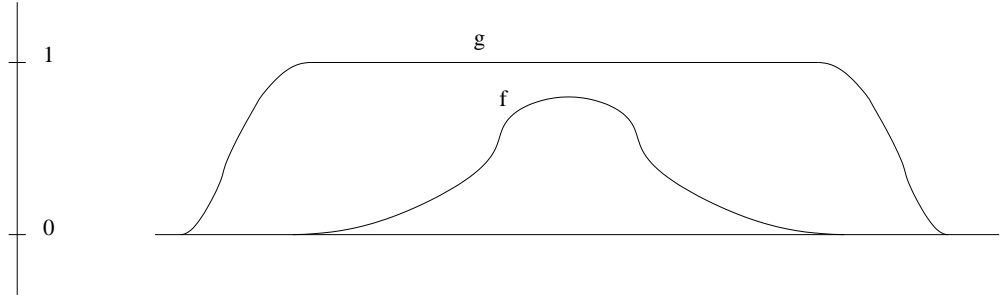}
\caption{\label{bump} A  bump function  $f$ of the Thom class representative and a characteristic function $g$}
\end{figure}

%Thus we can use $j_*\circ z_m^*$ to construct a section on forms.
To construct a section on the level of forms, we first notice that
the Thom class can be represented by using a bump function $f$ so
that  if $X^{m_i}$ is given locally  on $U$ by the equations $x_k =
\dots =x_{N}=0$

\begin{equation}
\T(1)|_U= f dx_k \wedge \dots \wedge dx_N
\end{equation}
where $f$ is a bump function along the fiber that can be chosen such
that $supp(f)$, the support of $f$, lies strictly inside the tubular
neighborhood and moreover $supp(f)$ lies strictly inside this
neighborhood. We consider a ``characteristic function'' $g$ of an
open subset $U$ with $supp(f) \subset U \subset Tub(N)$ inside the
tubular neighborhood, see Figure \ref{bump}. Notice that
$fg(x)=f(x)$. We let $\g$ be a 0--form with compact vertical support
whose restriction to the fiber is given by $g$.

For any form $\omega\in \Omega^*(X)$, we define
\begin{equation}
i_{ms}(\omega):= j_*(\g \pi_m^*(\omega))
\end{equation}

Then
\begin{equation}
i_m^*(j_*(\g \pi_m^*(\omega)))=z^*_m (j^*(j_*(\g \pi_m^*(\omega))))
=z^*_m(\g)z^*_m( \pi^*_m(\omega))= \omega +d\tau
\end{equation}

\begin{rmk}
Actually $i_*(\omega):= j_*(\T(w))$ is divisible by $j_*(\T(1))$:
locally on a coordinate neighborhood $U$.
\begin{eqnarray}
\label{dividebythomeq}
i_{m*}(\omega)&=&j_*(\T(w))=j_*(\pi^*_m(\omega)\wedge\Theta)|_U\nn\\
&=& f\pi^*(\omega)|_U \wedge dx_k \wedge \dots \wedge dx_N\nn\\
&=&fg\pi^*(\omega)|_U \wedge dx_k \wedge \dots \wedge dx_N\nn\\
&=&i_{ms}(\omega)|_U \wedge \Theta|_U
\end{eqnarray}
\end{rmk}

\begin{thm}
With $i_{ms}$ and $i_{m*}$ as defined above the following equation
holds on the level of forms.
\begin{eqnarray}
\label{formequation} \omega_{m_1}*\omega_{m_2} &:=& \check
e_{m_3*}(e^*_1(\omega_{m_1})e_2^*(\omega_{m_2})\Form(\e(\R(\bm))))\nn\\
&=& \coeff{t^r}\left\{ \right. \check
\imath_3^*[i_{s1}(\omega_{m_1})i_{s2}(\omega_{m_2})
i_{s1}(\Form(\et(S_{m_1})))i_{s2}(\Form(\et(S_{m_1}))) \nn\\
&&\quad \left. i_{s3}(\Form(\et(S_{m_3}))\Form(\et(\ominus
N_{X^{\bm}/X}))\Form(\e(N_{X^{\bm}/X^{m_3}})))]\right\}+d\tau
\end{eqnarray}
for some exact form $d\tau$, where $\Form(v)$ is a closed form
representative of the class $v$.
\end{thm}

\begin{proof}
Completely parallel to the proof of Proposition \ref{mainprop},
since we have established all equalities up to homotopy, that is up
to exact forms.
\end{proof}
\begin{cor}
The three point functions coincide with the ones induced by
(\ref{proddefeq}). That is if $\Form$ denotes the lift of a class to
a form and $\Form(v_{m_i})=\omega_{m_i}$ then

\begin{eqnarray}
\langle \omega_{m_1}*\omega_{m_2},\omega_{m_3}\rangle&:=&\int_X \omega_{m_1}*\omega_{m_2} \wedge \omega_{m_3}\nn\\
&=&
\int_X \Form (v_{m_1}*v_{m_2}) \wedge \omega_{m_3}\nn\\
&=& (v_{m_1}*v_{m_2}\cup v_{m_3})\cap [X]\\
&=&\langle v_{m_1}*v_{m_2},v_{m_3}\rangle
\end{eqnarray}
where $[X]$ is the fundamental class of $X$ and hence the three
point functions are independent of the lift.
\end{cor}
\begin{proof}
Straightforward by Stokes.
\end{proof}

\section{Formal and non--formal calculations}
\label{formalsection}

 In this section, we present some formal and
non--formal calculations. This will allow us to make contact with
the formal argument of Chen-Hu \cite{ChenHu} who used fractional
Thom forms in their arguments in establishing a deRham model for the
Chen--Ruan cohomology of Abelian quotients. We have been informed by
H.H.~Tseng \cite{Tseng} that he is working on making the formal part
of arguments rigorous in a ${\mathbb C}^*$ equivariant setting, a
result which would be great to have.

The ultimate aim would be to rigorously establish the following
presentation of the product for the various functors $\F$ without
recourse to the deRham theory or sections.

\begin{equation}
\label{maineq} \check \imath^{3*}(i_{1*}(v_{m_1
}\sigma_1)i_{2*}(v_{m_2}\sigma_2)\check
\imath_{3*}(\tilde\sigma_3))= \check
e_{3*}(e_1^*(v_{m_1})e_2^*(v_{m_2})\e(\R(\bm)))
\end{equation}

This is: can one answer the following questions?

\begin{qu}
Can one find elements $\sigma_i$, $\tilde \sigma_i$ such that the
equation (\ref{maineq}) holds?
\end{qu}

\begin{qu}
Is there a setting in which the  $\sigma_1,\sigma_2,\tilde \sigma_3$
form a co--cycle or better even a trivial co--cycle?
\end{qu}

In a formal sense this can be done as we show below, but we are
still lacking a rigorous setting. Of course the preceding paragraphs
do give rigorous results using the existence of sections for the
deRham and the cyclic setting.

We first notice that for $r=\rk(\R(\bm))$ the r.h.s. of
(\ref{maineq}) is
\begin{equation}
\evalf{r}[e_{3*}(e_1^*(v_1)e_2^*(v_2)\et(\R(\bm)))]
\end{equation}

\subsection{A rigorous excess intersection calculation}
 We calculate in $\F(X,G)[[t]]$. Using the excess intersection
formula \cite{FL,Quillen} on
\begin{equation}
\begin{matrix}
X^{\bm}&\stackrel{\check e_3}{\longrightarrow}&X^{m_3^{-1}}\\
\hskip-4cm (e_1,e_2,\check e_3)\circ (\Delta, id)\circ
\Delta\downarrow &&
\hskip2cm\downarrow (\check \imath_3,\check \imath_3,\check \imath_3)\circ (\Delta, id)\circ \Delta\\
X^{m_1}\times X^{m_2}\times X^{m^{-1}_3} &\stackrel{(i_1,i_2,\check
\imath_3)}\longrightarrow&X\times X\times X
\end{matrix}
\end{equation}
which has excess bundle $N_{X^{m_1}/X}\oplus N_{X^{m_2}/X} \oplus
N_{X^{m^{-1}_3}/X}|X^{\bm}\ominus N_{X^{\bm}/X^{m^{-1}_3}}$ we can
transform the l.h.s.\ of equation (\ref{maineq}) as follows:
\begin{eqnarray}
\label{lhseq} l.h.s. (\ref{maineq})&=&
\check \imath^*_3[i_{1*}(v_{m_1 }\sigma_1)i_{2*}(v_{m_2}\sigma_2)\check \imath_{3*}(\tilde \sigma_3)]\nn\\
&=& e_{3*}[e_1^*(v_{m_1 }\sigma_1\e(N_{X^{m_1}/X}))e_2^*(v_{m_2
}\sigma_2\e(N_{X^{m_2}/X}))\nn\\
&&\quad e_3^*(\tilde \sigma_3\e(N_{X^{m_3^{-1}/X}}))
\e(\ominus N_{X^{\bm}/X^{m_3^{-1}}})]\nn\\
 &=&\evalf{k}\left\{
\check e_{3*}[e_1^*(v_{m_1 }\sigma_1\et(N_{X^{m_1}/X}))e_2^*(v_{m_2
}\sigma_2\et(N_{X^{m_2}/X}))\right.\nn\\
&&\left.\quad \check e_3^*(\tilde
\sigma_3\et(N_{X^{m^{-1}_3}/X}))\et(\ominus
N_{X^{X^{\bm}/m_3^{-1}}})]\right\}
\nn\\
\end{eqnarray}
 $k=\rk(N_{X^{m_1}/X}\oplus N_{X^{m_2}/X} \oplus
N_{X^{m_3}/X}|X^{\bm}\ominus N_{X^{\bm}/X^{m_3^{-1}}})$.

While the r.h.s.\ can be transformed to
\begin{eqnarray}
\label{rhseq} r.h.s.(\ref{maineq})&=&\evalf{r}\left\{\check
e_{3*}[e_1^*(v_{m_1}\et(S_{m_1}))
e_2^*(v_{m_2}\et(S_{m_2}))\right.\nn\\
&&\left. \quad e_3^*(v_{m_3}\et(S_{m_3}))\et(\ominus
N_{X^{\bm}/X})]\right\}\nn\\
&=&\evalf{r} \left\{ \check
e_{3*}[e_1^*(v_{m_1}\et(S_{m_1})\et(N_{X^{m_1}/X})\et(\ominus
N_{X^{m_1}/X}))
\right.\nn\\
&&\left.\quad e_2^*(v_{m_2}\et(S_{m_2})\et(N_{X^{m_2}/X})\et(\ominus
N_{X^{m_2}/X}))\right.\nn\\
&&\left. \quad e_3^*(\et(S_{m_3}))e_3^*(\et (\ominus
N_{X^{m_3}/X}))\et(\ominus
N_{X^{\bm}/X^{m^{-1}_3}})]\right\}\nn\\
&=&\evalf{r} \left\{\check e_{3*}[e_1^*(v_{m_1}\et(\ominus
S_{m^{-1}_1})\et(N_{X^{m_1}/X}))
\right.\nn\\
&&\left.\quad e_2^*(v_{m_2}\et(\ominus S_{m_2})\et(N_{X^{m_2}/X}))
\check e_3^*(\et(\ominus S_{m^{-1}_3}))\et(\ominus
N_{X^{\bm}/X^{m^{-1}_3}})]\right\}\nn\\
\end{eqnarray}
$r=\rk(\R(m))$.

\subsection{A formal solution using fractional Euler classes}
Comparing the two sides that is equations (\ref{lhseq}) and
(\ref{rhseq}) one is tempted to set:
\begin{eqnarray}
\sigma_{1,t}&=&\et(\ominus S_{m^{-1}_1})=\et(S_{m_1})\et(\ominus
N_{X^{m_1}/X})\nn\\
\sigma_{2,t}&=&\et(\ominus S_{m^{-1}_2})=\et(S_{m_2})\et(\ominus
N_{X^{m_2}/X})\nn\\
\tilde \sigma_{3,t}&=&\et(\ominus S_{m^{-1}_3}\ominus
N_{X^{m_3^{-1}}/X})=\vee^*(\et(S_{m_3})\et(\ominus N_{X^{m_3}/X})^2)
\end{eqnarray}
and then use a kind of evaluation map that is set
$\sigma_i=\evalf{\vr(\sigma_i)}(\sigma_{i,t})$ and
 $\tilde \sigma_3:= \evalf{\vr(\tilde \sigma_3)}(\sigma_{3,t})$
where $\vr$ denotes the virtual rank. This is, however, not
possible, since it is not clear that the respective power series
converges for $-1$ nor is it clear what the coefficient at a
rational power or a negative virtual rank means.

\subsection{Adjoining formal symbols}
 Let $\fS$ be a collection of elements of rational
$K$--theory $K_{\mathbb Q}(Y)$.

We will think of the formulas first in $\F(Y)[\fS]$ and write
elements $\fS$ by using the formal symbols $\fe(x)$ (one should
think ``$\fe(x)=\evalf{\vr(x)}(x)$'')

We can see that we can ``solve'' the equation (\ref{maineq}) if we
{\em formally} set
\begin{eqnarray}
\label{eqnsols} \sigma_i&=&\fe(S_{m_i})\fe(\ominus
N_{X^{m_i}/X})%=\fe(\ominus S_{m_i^{-1}})
\nn\\
\tilde \sigma_3&=&\vee^*(\fe(S_{m_3})\fe(\ominus N_{X^{m_3}/X})^2)
%=\fe(\ominus S_{m_3}\ominus N_{X^{m_3}/X})
\\
\end{eqnarray}
as we explain in the following.

 One would like to add certain relations of the form
\begin{enumerate}
\item Enlarging $\fS$ to the semi--group it generates
\begin{equation}
\fe(x)\fe(y)-\fe(x\oplus y)
\end{equation}

\item If $x+y=E$ with $E$ a bundle
\begin{equation}
\fe(x)\fe(y)-\e(E)
\end{equation}

\end{enumerate}
We denote by $\FF{Y}{\fS}$ the ring obtained by modding out by the
relations above.

But one has to be careful with negative bundles, i.e.\ $\fe(\ominus
E)$, since these will be morally the inverses to nilpotent elements
and hence if we were to localize, we would obtain the zero ring.

Looking at our equations, we would like to have $Y=I(X)$ with the
maps $i_m^*$ $\fS=\{S_m,\ominus N_{X^m/X}|m\in G\}$ but then using
the relations above, we would get into trouble with
$\fe(S_m)\fe(S_{m^{-1}})\fe(\ominus N_{X^m/X})$.

What we will formally do is to view $\fe(\ominus N_{X^m/X})$ as
division, when it is possible, as we demonstrate below.

On the other hand, there is no problem adjoining only the $S_m$.

Therefore we will consider adjoining two sets of variables
$\fS_1:=\{S_m\}$ and $\fS_2:=\{\ominus N_{X^m/X}\}$. Then we will
consider the formulas to live in  $\FF{Y}{\fS_1}[{\mathfrak S}_2]$,
where we think of $\fS_2$ as formal division operators when defined
and also use the convention:

\begin{enumerate}
\item Pull--back: if $i^*(x)$ and $x$ in $K_{\mathbb Q}(Y)$ for some
morphism $i$: $i^*(\fe(x))=\fe(i^*(x))$
\end{enumerate}

\subsection{Formal manipulations: divisions give formal sections}
For a given inclusion $i:Y\to X$, we sometimes can construct sections
$i_s$ of $i^*$. Notice that $i_*$ is not quite a section due to the
self intersection formula:

\begin{equation}
i^*(i_*(a))=a\e(N_{X/Y})
\end{equation}
this is why we formally set
\begin{equation}
``i_s(a):= i_*(a\fe(\ominus N_{X/Y})) \text{''}
\end{equation}

Indeed, then
\begin{equation}
``i^*(i_s(a)):=i^*(i_*(a\fe(\ominus
N_{X/Y})))=a\e(N_{X/Y})\fe(\ominus N_{X/Y})=a\text{''}
\end{equation}

Notice that if $i_s$ is indeed a section:
\begin{equation}
i_*(ab)= i_*(i^*(i_s(a))b)=i_s(a)i_*(b)
\end{equation}
and hence
\begin{equation}
i_*(a)=i_s(a)i_*(1)
\end{equation}

So that we see that if there are sections:
\begin{equation}
i_s(a)=i_*(a)/i_*(1)
\end{equation}
and hence indeed the division operation is well justified and the
formal calculation is valid. This was the case in
\S\ref{cyclicsection} and \S\ref{deRhamsection}, see in particular
equation (\ref{dividebythomeq}).

In the notation above, the l.h.s.\ of (\ref{maineq}) after
substitution of (\ref{eqnsols}) becomes
\begin{equation}
\label{formallhseq} i_3^*[i_{s1}(v_{m_1})i_{s2}(v_{m_2})i_{s1}
(\fe(S_{m_1}))i_{s2}(\fe(S_{m_1}))i_{s3}(\fe(S_{m_3})\fe(\ominus
N_{X^{m_3}/X}))]
\end{equation}
while the r.h.s.\ of (\ref{maineq})can be written as

\begin{equation}
\label{formalrhseq}
 \check
e_{3*}(e^*_1(v_{m_1})e^*_2(v_{m_2})e^*_1(\fe(S_{m_1}))
e^*_2(\fe(S_{m_2}))e^*_3(\fe(S_{m_3}))\fe(\ominus N_{X^{\bm}/X}))
\end{equation}

\subsubsection{Cocycles} Notice since $S_m\oplus
\vee^*(S_{m^{-1}})=N_{X^m/X}$ we formally have that
\begin{multline}
``\check \imath_{3s}
(\fe(S_{m^{-1}_3})i_{3s}(\fe(S_{m_3})\fe(\ominus N_{X^{m_3}/X})))
\\=i_{s3}(\e(N_{X^{m_3}/X})\fe(\ominus N_{X^{m_3}/X}))=1\text{''}
\end{multline}
This if we set $s(m):=i_{ms}(\fe(S_m))$ and let $\gamma:=ds$, that
is $\gamma(m_1,m_2) = s(m_1)s(m_2)/s(m_1m_2)$, then the l.h.s.\ of
(\ref{maineq}) becomes
\begin{equation}
i_3^*[i_{s1}(v_{m_1})i_{s2}(v_{m_2})\gamma(m_1,m_2)]
\end{equation}
with a trivial co-cycle
\begin{equation}
``\gamma(m_1,m_2)=i_{1s}(\fe(S_{m_1}))i_{2s}(\fe(S_{m_2})) \check
\imath _{3s}(\fe(S_{m^{-1}_3}))^{-1}\text{''}
\end{equation}

 This formal equation is very important,
since it makes contact with the algebraic problem posed and studied
in \cite{orb} called the re--construction problem in
\cite{wisc,orb,hilb}. This program has previously been very useful
for symmetric products \cite{sq} and singularities with symmetries
\cite{orbsing}.

\subsection{Positive fractional Euler-classes}
\label{posEulerfrakpar} Unlike the negative fractional Euler classes
$\fe(\ominus N)$, we can make the Euler classes of positive rational
combinations of bundles rigorous.

First we notice that by the splitting principle \cite{Hirz,FL}, we
can make a ring extension in which all the constituent bundles
split. Then we are left with classes of the form $\fe(\frac{k}{r}
\L)$ that is fractional line bundles. Let $1+u=\e(\L)$ then we can
easily adjoin r--th roots to the extension of $R'$ of $R:=\F(X,G)$
in which all isotypical components of the $N_{X^{m}/X}$ split by
passing to $R'[w]/(w^r-u)$. After adjoining all $|m|$--th roots of
the various $\L_{m,k,i}$, where the $\L_{m,k,i}$ are the bundles
that split $W_{m,k}$, we can simply set
\begin{equation}
\label{fractionaldefeq}
 \fe(S_m):=\prod_{k\neq 0,i} w_{m,k,i}^{k/|m|}
\end{equation}
In this large ring $R$ is a subring and hence we can read of
formulas on this subring analogously to the procedure used in the
splitting principle.

\subsection{Admissible functors}
Here we collect the formal properties of the functors $\F$ we used
in our formal calculations.

\begin{df}
Let $\F$ be a functor together with an Euler-class  $\et$ which has
the following properties
\begin{enumerate}
\item $\F$ The Euler class $\et$ is defined for elements of rational K-theory and is multiplicative and takes values in $\F(X)[[t]]$.
\item $\F$ is contravariant,
i.e.\ it has pullbacks and the Euler-class is natural with respect to these.
\item $\F$ has push-forwards $i_*$ for closed embeddings $i:X\hookrightarrow Y$.
%!!closed???
\item $\F$ has an excess intersection formula for closed embeddings.
That is we have an evaluation morphism
$\e:=\evalf{r}(\eut{\F}):\F(X)[[t]]\to \F(X)$ such that for the
Cartesian squares
\begin{equation}
\begin{matrix}
Z&\stackrel{e_2}{\rightarrow}&Y_2\\
\downarrow e_1 & & \downarrow i_2\\
Y_1&\stackrel{i_1}\rightarrow&X
\end{matrix}
\end{equation}
we have the following formula
\begin{equation}
i_2^*(i_{1*}(a))=e_{2*}(e^*_1(a)\eps_j)
\end{equation}
where $\eps:=\e(E)$ with $E$ the excess bundle $E:=N_{Y_1/X}|_{Z}
\ominus N_{Z/Y_2}$ and $r$ is its rank.
\end{enumerate}

We call such a functor {\em admissible}.
\end{df}
All the functors $\F$ studied above are admissible and the
calculations of this section ---formal and non--formal--- carry over
to admissible functors. Actually deRham forms are admissible up to
homotopy, see below, so that {\em mutatis mutandis} we can use the
same formal arguments on the level of forms.

\subsection{Forms as an admissible functor, fractional Thom classes}

In this case, we have an Euler class and all the properties of an
admissible functor are valid on the chain level - up to homotopy,
that is up to closed forms.

\begin{enumerate}
\item The Thom push--forward on the chain level induces the
push--forward in cohomology induced by the Poincar\'e pairing, since
the Thom class and the Poincar\'e dual can be represented by the
same form \cite{BottTu}.
\item The projection formula holds, since the pull--back of the Thom
class is the Euler class of the normal bundle \cite{BottTu}.
\item The excess intersection formula holds up to homotopy. Since it
holds in cobordism theory and cohomology \cite{Quillen} we know that
for closed $\omega$ the two forms $i_2^*i_{1*}(\omega)$ and
$e_{2*}e^*_1(\omega \Form(\e(E)))$ differ by a closed form.
\item
In particular, we can use the Thom pushforward and then the
divisibility of the push--forward by the Thom class gives us
sections.
\end{enumerate}
Hence we can make the same formal calculations as above. Notice that
since we indeed have sections as explained in \S\ref{deRhamsection},
we can avoid mention of $\fe(\ominus N_{X^m/X})$ and only have to
deal with $\fe(S_m)$. In particular equation (\ref{dividebythomeq})
shows that we can indeed divide $i_{m*}(a)$ by the Thom form
$i_{m*}(1)=\Theta_m$, which is how we defined $i_{ms}$.

Now using the formalism of \S\ref{posEulerfrakpar} and passing to a
local trivializing neighborhood $U$, where the line bundles
$\L_{m,k}$ have first Chern class represented by the forms
$dx_l,\dots dx_N$, then we get a form representative of $\fe(S_m)$

\begin{equation}
\Form(\fe(S_m))|_U=f^{k/|m|}r\prod_{k\neq0,i} (dx)^{k/|m|}
\end{equation}
which is the expression for the fractional Thom form that was used
in \cite{ChenHu} in their study of Abelian quotients.

What we have now is the generalization to an arbitrary group as well
as a trivialization of the co--cycles in terms of roots, thus
completing that (re)--construction program of \cite{wisc,orb} in the
deRham setting of global quotients. The surprising answer is that
there is always a stringy multiplication arising from a co--cycle
that is trivializable in a ring extension obtained by adjoining
roots. In particular the formulas (\ref{formallhseq}) and
(\ref{formalrhseq}) can be made sense of and the  co--cycle that
appears yields the stringy orbifold product.

%something about fractional Thom classes.

\end{document}